\documentclass[a4paper]{amsart}

\usepackage{amssymb}
\usepackage{stmaryrd}

\def\frak{\mathfrak}
\def\Bbb{\mathbb}
\def\Cal{\mathcal}
\def\sideremark#1{\ifvmode\leavevmode\fi\vadjust{\vbox to0pt{\vss
 \hbox to 0pt{\hskip\hsize\hskip1em
 \vbox{\hsize3cm\tiny\raggedright\pretolerance10000
 \noindent #1\hfill}\hss}\vbox to8pt{\vfil}\vss}}}%

\swapnumbers
\numberwithin{equation}{subsection}

\newtheorem{prop}[subsection]{Proposition}
\newtheorem*{prop*}{Proposition}
\newtheorem{thm}[subsection]{Theorem}
\newtheorem*{thm*}{Theorem}
\newtheorem{lem}[subsection]{Lemma}
\newtheorem*{lem*}{Lemma}

\newtheorem*{kor*}{Corollary}
\theoremstyle{remark}

\newtheorem*{remark*}{Remark}
\newtheorem{example}[subsection]{Example}
\newtheorem*{example*}{Example}

\newcommand{\ad}{\operatorname{ad}}
\newcommand{\Ad}{\operatorname{Ad}}
\newcommand{\Adb}{\operatorname{\underline{Ad}}}
\renewcommand{\exp}{\operatorname{exp}}

\newcommand{\Fl}{\operatorname{Fl}}
\newcommand{\x}{\times}
\renewcommand{\o}{\circ}
\newcommand{\ddt}{\tfrac{d}{dt}|_{t=0}}
\newcommand\assoc[2]{\llbracket #1,#2\rrbracket}
\newcommand{\dev}{\operatorname{dev}}
\newcommand{\fg}{\frak g}

\let\ccdot\cdot
\def\cdot{\hbox to 2.5pt{\hss$\ccdot$\hss}}

\newcommand{\al}{\alpha}

\newcommand{\de}{\delta}

\newcommand{\la}{\lambda}
\newcommand{\rh}{\rho}
\newcommand{\om}{\omega}
\renewcommand{\phi}{\varphi}
\newcommand{\ph}{\varphi}
\newcommand{\ps}{\psi}

\newcommand{\ze}{\zeta}

\newcommand{\Om}{\Omega}

\def\R{\Bbb R}
\def\>{\rightarrow}
\def\g{\frak g}
\def\p{\frak p}

\def\C{{\Cal C}}

\def\S{{\Cal S}}
\def\ddt#1{\left.\tfrac{d}{dt}\right\vert_{#1}}
\def\smatrix#1{\left(\smallmatrix #1 \endsmallmatrix\right)}
\def\c#1{c^{b_#1,X_#1}}
\def\d#1{^{(#1)}}

\begin{document}
\title{On distinguished curves in Parabolic Geometries} 
\author{Andreas \v Cap, Jan Slov\'ak, and Vojt\v ech \v Z\'adn\'\i k}
\date{August 6, 2003}

\address{A.C.: Institut f\"ur Mathematik, Universit\"at Wien,
Strudlhofgasse~4, A--1090 Wien, Austria, and International Erwin
Schr\"odinger Institute for Mathematical Physics, Boltzmanngasse 9, A-1090
Wien, Austria\newline\indent J.S. and V.\v Z.: 
Department of Algebra and Geometry,
Masaryk University, Jan\'a\v ckovo n\'am. 2a, 662~95~Brno, Czech Republic}

\email{andreas.cap@esi.ac.at, slovak@math.muni.cz, zadnik@math.muni.cz}

\begin{abstract}
  All parabolic geometries, i.e.~Cartan geometries with homogeneous
  model a real generalized flag manifold, admit highly interesting
  classes of distinguished curves. The geodesics of a projective class
  of connections on a manifold, conformal circles on conformal
  Riemannian manifolds, and Chern--Moser chains on CR--manifolds of
  hypersurface type are typical examples. We show that such
  distinguished curves are always determined by a finite jet in one
  point, and study the properties of such jets. We also discuss the
  question when distinguished curves agree up to reparametrization and
  discuss the distinguished parametrizations in this case. We give a
  complete description of all distinguished curves for some examples of
  parabolic geometries.
\end{abstract}

\maketitle

Elie Cartan's idea of `generalized spaces' as curved analogs of Felix
Klein's geometries (i.e. homogeneous spaces) is a well understood
geometrical concept, which, for a Lie subgroup $P\subset G$,
generalizes the Maurer--Cartan form on the total space of the
principal $P$--bundle $G\to G/P$ to Cartan connections on principal
$P$--bundles, see e.g.~the introductory book \cite{Sha}. The concept
of \textit{parabolic geometries} refers to those cases where $P$ is a
parabolic subgroup in a (real or complex) semisimple Lie group $G$.
In \cite{Fef}, Ch. Fefferman initiated a program to exploit the
representation theory of parabolic subgroups in semisimple Lie groups
in order to understand invariants of geometric structures like
CR--geometries, projective geometries, or conformal Riemannian
geometries. This approach has proved to be extremely powerful. First
of all, all parabolic geometries can be described in terms of weaker
analogies of classical G--structures on smooth manifolds and,
similarly to the examples mentioned above, all such structures give
rise to canonical normal Cartan connections, cf.  \cite{Tan, Mor, CS}. In
fact, these constructions express Cartan's method of equivalence
using the language of the modern representation theory and natural
cohomological reasoning. The existence of the Cartan connection
provides an effective calculus to deal with invariant objects, see
e.g.~\cite{CSS4} and the references therein. To large extent, the
understanding of the general (curved) geometries can be reduced to
properties of the homogeneous model, and thus to purely algebraic
questions.

The goal of this paper is to use this approach in order to understand
invariantly defined systems of distinguished curves for parabolic
geometries, which we call \textit{(generalized) geodesics}.  After
recalling basic concepts of parabolic geometries, geodesics are
introduced and discussed along the lines of the classical approach in
affine geometry, which uses the development of curves. This approach
may be found in similar context already in \cite{Sha} and in
\cite{Koch3}.  In this way, many aspects of the study of the curves
are reduced to the case of the homogeneous model. Thus the original
`smooth' question on curved manifolds can be transformed to an
`algebraic' problem, which is discussed in Section \ref{2}. In
particular, we obtain estimates on the order of jets necessary to
determine a geodesic, and this approach 
also leads to an algebraic description of
all jets of geodesics in a point. The third section is devoted to the
study of possible reparametrizations in the class of geodesics.
Specializing the general results to $|1|$--graded Lie algebras, we
obtain generalizations of some well known results on conformal,
projective, and quaternionic geometries (see e.g.  \cite{BE}).  The
final section provides further refinements for specific classes of
curves, see in particular Theorems \ref{4.2} and \ref{4.3}.

\medskip 
\noindent{\bf Acknowledgments.}    
Part of the work was done during a stay of the second author at the
University of Adelaide under an ARC financial support, and his
discussions with Michael Eastwood were most helpful and illuminating.
First author supported by project P15747 of the FWF. The second and
third authors acknowledge the support from GACR, Grant Nr.
201/02/1390.

\section{General concepts}\label{1}

\subsection{Parabolic geometries}\label{1.1} 
Let us briefly recall the basic facts, more details can be found in
\cite{CS-weyl} or \cite{Sha}, and the references therein.  Let $G$ be
a real semisimple Lie group with Lie algebra $\frak g$, and $P\subset
G$ a parabolic subgroup with Lie algebra $\frak p$. A (real)
\textit{parabolic geometry} $(\Cal G, \om)$ of type $(G,P)$ is a
principal bundle $\Cal G$ with structure group $P$ over a manifold
$M$, equipped with a smooth one--form $\om\in\Om^1(\Cal G,\frak g)$,
which satisfies

(1) $\om(\ze_Z)(u)=Z$ for all $u\in \Cal G$ and fundamental fields $\ze_Z$,
$Z\in\frak p\subset \frak g$, i.e.~$\om$ reproduces the generators of
fundamental vector fields,

(2) $(r^b)^*\om= \Ad(b^{-1})\o \om$ for all $b\in P$, i.e. $\om$ is
$P$--equivariant with respect to the adjoint representation, and

(3) $\om|_{T_u\Cal G}: T_u\Cal G\to \frak g$ is a linear isomorphism for
all $u\in \Cal G$, i.e. $\om$ is an absolute parallelism on $\Cal G$.

The curvature of a parabolic geometry $(\Cal G, \om)$ is the horizontal
two--form $K\in\Om^2(\Cal G, \frak g)$ defined by the structure equations
$$
K = d\om +\tfrac12[\om,\om]\text{, i.e.~}
K(\xi,\eta)=d\om(\xi,\eta)+[\om(\xi),\om(\eta)].
$$

Clearly, the Maurer--Cartan form $\om$ on the principal fiber bundle
$G\to G/P$ is a parabolic geometry and the structure equations say
that this geometry is \textit{flat}, i.e.~its curvature vanishes
identically. $(G\to G/P,\om)$ is called the \textit{homogeneous model}
for parabolic geometries of type $(G,P)$. 

Morphisms between Cartan geometries $(\Cal G,\om)$ and $(\Cal
G',\om')$ are principal fiber bundle morphisms $\ph: \Cal G\to \Cal
G'$ such that $\ph^*\om'=\om$. It is quite elementary to prove that a
geometry is locally isomorphic to its homogeneous model if and only if
its curvature vanishes identically, cf. \cite{Sha}.

Each smooth (left) action of the structure group $P$ on a smooth
manifold $S$ leads to a functor $\Cal S$ on the category of Cartan
geometries of type $(G,P)$. The value of $\Cal S$ on $(\Cal G,\om)$ is
the associated fiber bundle $\Cal G\x_P S$ with respect to the action
of $P$ while a morphisms $\ph:(\Cal G,\om)\to (\Cal G',\om')$ induces
the fiber bundle morphism $\ph\x_P \operatorname{id_S} : \Cal G\x_P
S\to \Cal G'\x_P S'$. We call these bundles \textit{natural bundles}.
Moreover, this construction is functorial in the smooth action entry,
because each equivariant mapping $\al:S\to S'$ induces the fiber
bundle mapping $\operatorname{id_{\Cal G}}\x_P \al : \Cal G\x_P S\to
\Cal G\x_P S'$. Thus we have got a bifunctor on Cartan geometries and
smooth left actions with values in fiber bundles.

In particular, linear representations of $P$ lead to functors valued in
vector bundles and their linear morphisms and the bifunctoriality of the
construction extends all natural constructions like pairings,
decompositions, and tensor products of representations to the natural bundles.
Of course, all this is the obvious restriction of the usual functorial
constructions over all principal fiber bundles to the category of Cartan
geometries.

A central example, which also illustrates the role of the Cartan
connection, is given by the representation of $P$ on $\frak g/\frak p$
induced by the adjoint representation. This leads to the functor $\Cal
G\x_P \frak g/\frak p$, and via the Cartan connection $\om$ this
associated bundle can be identified with the tangent bundle $TM$.
Indeed, since $\om$ defines an absolute parallelism, there are the
corresponding `constant' vector fields $\om^{-1}(X)\in \Cal X(\Cal G)$
for all $X\in \frak g$, defined by $\om(\om^{-1}(X)(u))=X$ for all
$u\in \Cal G$. Denoting by $\assoc{u}{X+\frak p}$ the class in $\Cal
G\x_P\frak g/\frak p$ of $(u,X+\frak p)\in\Cal G\x\frak g/\frak p$ and
by $\pi:\Cal G\to M$ the bundle projection, one immediately verifies
that $\assoc{u}{X + \frak p}\mapsto T\pi(\om^{-1}(X)(u))$ defines the
claimed isomorphism.

For any parabolic subalgebra $\frak p\subset\frak g$, there is a
grading $\fg_{-k}\oplus\dots \oplus\fg_k$ of $\fg$ such that $\frak
p=\fg_0\oplus\dots\oplus\fg_k$, and $\frak
p_+=\fg_1\oplus\dots\oplus\fg_k$ is the nilradical of $\fg$, see
\cite{Yam, CS}. In particular, this implies that $\fg_0$ is a
reductive Levi component for $\frak p$.  Hence we obtain an
identification $\frak n=\fg_{-k}\oplus\dots\oplus\fg_{-1}$ with $\frak
g/\frak p$, which is an isomorphism of $P$--modules if we endow $\frak
n$ with the `truncated' adjoint action $\Adb$. Via the
Killing form, one further obtains an identification of $\frak n^*$ with
$\frak p_+$, which induces the identification of the cotangent bundle
$T^*M$ with $\Cal G\x_P \frak n^*$. Thus all tensor bundles over $M$
are identified with the natural bundles coming from tensor products of
the representations $\frak n$ and $\frak n^*$. Moreover, the right
hand ends $\frak g^i=\frak g_i\oplus\dots\oplus \frak g_k$ define a
$P$--invariant filtration of $\fg$. Hence we obtain natural subbundles
$T^iM\subset TM$ for all $i<0$. The resulting filtration
$$
TM=T^{-k}M\supset T^{-k+1}M\supset \dots \supset T^{-1}M\supset 0
$$ 
is the most importing object underlying a parabolic geometry. This
filtration is trivial for $|1|$--graded algebras and we call such parabolic
geometries \textit{irreducible}. 

A very special case of the construction of natural bundles is the
choice $S=G$ with the left action of $P$ on $G$ given by the group
multiplication.  This leads to the principal fiber bundle $\tilde{\Cal
  G} = \Cal G\x_P G$ with the principal action given by the usual
right multiplication in $G$ and the canonical inclusion ${\Cal
  G}\subset \tilde{\Cal G}$, $u\mapsto \assoc{u}{e}$, where $e\in G$
is the unit element.  Now, the Cartan connection $\om$ extends
uniquely to a $G$--equivariant one--form $\tilde\om\in\Om^1(\tilde
{\Cal G},\frak g)$ reproducing the fundamental vector fields. 
One easily verifies that $\tilde\om$ is a
principal connection on $\tilde{\Cal G}$. Whenever we have a left
action of $P$ on some manifold $S$ which is the restriction of a left
action of $G$, then we may view the natural bundle $\Cal G\x_P S$ also
as $\tilde{\Cal G}\x_G S$. Hence on any natural bundle of this type,
there is a canonical connection induced by $\tilde\om$. Of course, if
we consider restrictions of $G$--representations to $P$, then the
resulting natural vector bundles, which are usually called
\textit{tractor bundles}, are equipped with canonical linear
connections. 

\subsection{Development of curves}\label{1.2}
The notion of the development of curves is related to a particular
instance of natural bundles associated to restrictions of $G$--actions
to $P$, namely the case of the canonical left action on $G/P$. The
resulting space $\Cal S=\Cal G\x_P G/P=\tilde{\Cal G}\x_G G/P$ is
called \textit{Cartan's space} over the underlying manifold $M$ of the
Cartan geometry in question. Of course, $\Cal S\to M$ is a fiber
bundle with typical fiber $G/P$, and from \ref{1.1} we know that 
the parabolic geometry induces a canonical connection on this fiber
bundle. 

Another remarkable fact about $\Cal S$ is that for the point $o=eP\in
G/P$, and a point $x\in M$, all points $u\in\Cal G$ with $\pi(u)=x$
lead to the same class $O(x)=\assoc uo\in \Cal G\x_P G/P$. Hence we
obtain a canonical smooth section $O$ of $\Cal S\to M$ for every
parabolic geometry $(\Cal G\to M,\om)$ of type $(G,P)$. Moreover, the
vertical tangent bundle $V\Cal S$ can be identified with the
associated bundle $\Cal G\x_P T(G/P)$. Since the base point $o\in G/P$
is a fix point for the action of $P$, we see that the restriction of
$V\Cal S$ to the image $O(M)$ of the canonical section is the
associated bundle $\Cal G\x_P T_o(G/P)$. Since $T_o(G/P)$ is
canonically isomorphic with $\frak g/\frak p$ and $\Cal G\x_P(\frak
g/\frak p)$ is naturally isomorphic to $TM$, we get a canonical
isomorphism $V\Cal S|_{O(x)}\cong TM$. Thus we may view the Cartan's
space $\Cal S$ as a nonlinear version of the tangent bundle in which
the geometry in question is encoded by means of the local parallel
transport of the induced connection.  This point of view goes back to
Cartan, and it was developed further in an abstract way in the second
half of the 20th century (see e.g. \cite{Kolar}).

This canonical parallel transport provides a straightforward
generalization of the classical concept of the development of curves.
By composing with $O$, a curve $c:I\to M$ with $I=(a,b)\subseteq\Bbb R$
may be also viewed as a parametrized curve in $\Cal S$. Fixing $t_0\in
I$ we find a neighborhood $J$ of $t_0$ in $I$ on which the parallel
transport along $c:I\to M$ is well defined. Given $s\in J$, we may
follow the curve $O\o c$ from $t_0$ to $s$ and then follow the
parallel transport backward for time $t_0-s$ to return to the fiber
over $t_0$. More formally, we define a smooth curve $\dev(c,t_0)$ from
an open neighborhood of $0$ in $\Bbb R$ to $\Cal S_{c(t_0)}$ by
$\dev(c,t_0)(s):=\tilde c_s(s)$, where $\tilde c_s$ is the parallel
curve in $\Cal S$ lying over $t\mapsto c(t_0+s-t)$ with the initial
point $O(c(s))$. This curve is called the \textit{development} of $c$
at $t_0$. For a point $u\in \Cal G$ over $c(t_0)$, there is a unique
curve $\bar c(t)$ in $G/P$ mapping $0\in\Bbb R$ to $o\in G/P$ such
that $\dev(c,t_0)(t)=\assoc{u}{\bar c(t)}$. Any other choice for the
point in $\Cal G$ has the form $u\cdot b$ for $b\in P$, and for that
choice the curve changes to $\ell_{b^{-1}}\o\bar c$. 

Hence we conclude that each choice of a $P$--invariant class $\Cal C$
of curves which map $0\in\Bbb R$ to $o\in G/P$ leads to a
distinguished class of curves on all manifolds endowed with a Cartan
geometry of type $(G,P)$. We say that a curve $c$ on $M$ is a
\textit{distinguished curve of type $\Cal C$} at a point $c(t_0)\in
M$, if for some (and thus any) point $u\in\Cal G$ the curve $\bar c$
constructed above lies in $\Cal C$. 

The natural choices for such sets $\Cal C$ of curves, of course come
from one--parameter subgroups in $G$: For a subset $A\subseteq \frak g$,
we can define a class $\Cal C_A$ as $\{ t\mapsto b\exp(tX)P:X\in A, b\in
P\}$. So we take the one--parametric subgroups with generators in $A$,
allow them to be shifted by left multiplications with elements of $P$,
and project the resulting curves to $G/P$. Of course, for $X\in \frak
p$ this always leads to the constant curve $o$, so we may assume
$A\cap \frak p=\emptyset$. On the other hand, if we want to have
curves in all directions in the class $\Cal C_A$, then we have to
assume that the restriction of the projection $\fg\to\fg/\frak p$ to
$A$ is surjective. The most obvious choice for $A$ which satisfies
this requirement is $A=\frak n$. It should be noted that for $X\in
\frak g\setminus \frak p$ the curve $t\mapsto b\exp(tX)P$ does not lie
in $\Cal C_{\frak n}$ in general. Following the case of affine
geometry and since we are mainly interested in having sets of
distinguished curves which are as small as possible, we shall always
assume $A\subseteq \frak n$ in the sequel.

The parabolic subgroup $P\subset G$ always has a canonical closed
subgroup $G_0$ which corresponds to the Lie subalgebra
$\fg_0\subset\frak p$. This group turns out to be reductive, and it
can be characterized as the subgroup of those elements in $G$, whose
adjoint action preserved the grading of $\fg$. In particular, the
subspace $\frak n$ is stable under the adjoint action of $G_0$. Now
for $b\in G_0$ and $X\in\frak n$, we of course have
$b\exp(tX)=\exp(t\Ad_b X)b$, and thus $b\exp(tX)P=\exp(t\Ad_b
X)P$. Thus it is natural to restrict attention to $G_0$--invariant
subsets $A\subseteq\frak n$, and the corresponding distinguished curves
are called \textit{(generalized) geodesics} of type $\Cal C_A$.  We
often do not mention the type if $A=\frak n$.

The generalized geodesics of type $\Cal C_A$ are easily
described explicitly by means of the constant vector fields
$\om^{-1}(X)$. Let us consider the projection $c(t)$ of the flow line
$\Fl^{\om^{-1}(X)}_t(u)\in \Cal G$ to the manifold $M$. From the
construction of the principal connection $\tilde \om$ on $\tilde{\Cal
  G}$ one immediately concludes that the horizontal vectors for
$\tilde\om$ in points $u\in \Cal G$ are $\om^{-1}(X)(u)-\ze_{X}(u)$
for all $X\in \frak n$. Thus, the curve $t\mapsto
\Fl^{\om^{-1}(X)}_t(u)\cdot \exp(-tX)$ must be the horizontal lift of
$c$ to $\tilde{\Cal G}$.  Now, the induced parallel transport of an
element $\assoc{u}{\exp tX}\in \Cal S$ along $c$ is given at time $s$
by $\assoc{\Fl^{\om^{-1}(X)}_s(u)}{\exp(t-s)X}$ and it reaches exactly
the point $O(c(t))$ in the canonical embedding of $M$ into $\Cal S$ at
time $s=t$. But this exactly means that for each $X\in\frak n$ the
curve $t\mapsto\assoc{u}{\exp tX}$ is the development of the
projection of the flow line through $u$ of the constant vector field
$\om^{-1}(X)\in\frak X(\Cal G)$. Since the allowed developments for
curves in $\Cal C_A$ have the form $t\mapsto\assoc{u}{\exp tX}$ for
$u\in \Cal G$ and $X\in A$, we have proved the first part of:

\begin{prop}\label{1.3} 
Let $(p:\Cal G\to M,\om)$ be a parabolic geometry of type $(G,P)$ and
let $A\subseteq\frak n$ be a $G_0$--invariant subset.  

\noindent
(1) The geodesics of type $\Cal C_A$ on $M$ are exactly the
projections of flow lines of the constant vector fields
$\om^{-1}(X)\in\frak X(\Cal G)$ with $X\in A$. 

\noindent
(2) Let $(p':\Cal G'\to M',\om')$ be another parabolic geometry of type
$(G,P)$, $\ph:\Cal G\to\Cal G'$ be a morphism of parabolic geometries
covering $\ph_0:M\to M'$, and $c:I\to M$ a smooth curve. Then $c$ is a
geodesic of type $\Cal C_A$ if and only if $\ph_0\o c:I\to M'$ is a
geodesic of type $\Cal C_A$.  
\end{prop}
\begin{proof}
  The curve $c$ in $M$ is a geodesics if an only if
  $c(t)=p\o\Fl^{\om^{-1}(X)}_t(u)$ for some $X\in A$ and $u\in \Cal
  G$. Since $\ph^*\om'=\om$, we get 
$$
p'\o\Fl^{{\om'}^{-1}(X)}_t(\ph(u)) = p'\o \ph\o\Fl^{\om^{-1}(X)}_t(u) =
\ph_0\o p\o \Fl^{\om^{-1}(X)}_t(u)
$$
and the claim follows. 
\end{proof}

\begin{remark*}
  (1) Our definition of geodesics and their general description is
  valid for arbitrary Cartan geometries. Though this is not a
  parabolic geometry, we may thus illustrate it in the case of affine
  connections on manifolds (i.e. $G$ is the affine group $\Bbb R^m\x
  GL(m,\Bbb R)$ and $P=GL(m,\Bbb R)$). Here the complement $\frak
  n=\Bbb R^m$ is $P$--invariant, and so any Cartan connection $\om$ on
  $\Cal G$ splits into the soldering form $\om_{\frak n}\in\Om^1(\Cal
  G,\Bbb R^m)$ and the principal connection form $\om_{\frak
    p}\in\Om^1(\Cal G,\frak p)$. Thus a Cartan geometry equips the
  underlying manifold $M$ with the linear frame bundle $(\Cal G,
  \om_{\frak n})$ and the principal connection $\om_{\frak p}$ on
  $\Cal G$. The projections of flow lines of the constant vector
  fields are exactly the geodesics of the linear connection on $TM$
  induced by $\om$.  Part (1) of the Proposition recovers the
  classical fact that the geodesics are those curves whose
  developments are straight lines in $\Bbb R^m=G/P$. On the other
  hand, if we choose $A=\frak g\setminus \frak p$, then more curves
  appear. For example, the following curves are projections of shifts
  of one--parametric subgroups in the affine group to the plane $\Bbb
  R^2$: $y=x\operatorname{log}x$ through $(1,0)$,
  $y=\operatorname{e}^x$ through $(0,1)$, $y=x^\al$ through $(1,1)$,
  cf. \cite{Eas-affine}.

\noindent
(2) Exactly as in the homogeneous case, each choice of $u\in\Cal G$
defines local coordinates around its projection $p(u)\in M$. Consider
the mapping $X\mapsto p(\Fl^{\om^{-1}(X)}_1(u))$, which is well
defined on some neighborhood $U\subseteq\frak n$ of $0$. Choosing $U$
sufficiently small, this becomes a diffeomorphism onto its image, thus
gives rise to local coordinates on $M$. These are called \textit{\em
normal coordinates} for the Cartan geometry in question. Of course,
in the setting of (1), we recover exactly the usual normal coordinates
for affine  connections on manifolds in this way. More information and
a characterization of the normal coordinates can be found in
\cite{CS-weyl}. 

We may rephrase our definition in terms of normal coordinates as
follows: The geodesics of type $\Cal C_A$ are those curves which are
linearly parametrized straight lines through the origin with
directions in $A\subseteq \frak n$ in some normal coordinates. Again,
this generalizes the standard facts on affine connections.
\end{remark*}

\begin{example}\label{1.5}
Let us mention four well known examples of distinguished curves in
parabolic geometries:

\noindent
(1) $G=SL(m+1,\Bbb R)$, $P$ is the stabilizer of a line in $\Bbb
R^{m+1}$. Normal parabolic geometries of type $(G,P)$ are classical
projective structures on $m$--dimensional manifolds. Generalized
geodesics (of type $\Cal C_{\frak n}$) are exactly the geodesics of
all connections in the projective class. They are determined by their
$2$--jet in one point as parametrized curves, but already determined
by their direction in one point as unparametrized curves.

\noindent
(2) $G=SL(m+1,\Bbb H)$, $P$ is the stabilizer of a quaternionic line.
This choice leads to almost quaternionic geometries (the complex
version of which is dealt with in \cite{BE}).  Again generalized
geodesics are determined by their $2$--jet in one point, but they form
more complicated systems of curves than in the projective case, cf.
\cite{BE}.

\noindent
(3) $G=O(p+1,q+1)$, $P$ is the stabilizer of a null line. This leads
to conformal pseudo Riemannian geometries of signature $(p,q)$. Here
the (generalized) geodesics are the well known conformal circles,
which owe their name to the fact that for the homogeneous model with
signature $(n,0)$ one obtains all circles on the sphere. For general
signatures, the geodesics in null directions, which behave similarly
to the projective case, form an interesting subclass.

\noindent
(4) $G=SU(p+1,q+1)$, $P$ the stabilizer of a (complex) null line. This
Hermitian analog of (3) leads to non--degenerate CR--structures of
hypersurface type with signature $(p,q)$. Here the Lie algebra is
$2$--graded and the geodesics of type $\Cal C_{\frak g_{-2}}$ are the
well known Chern--Moser chains. 
\end{example}

\section{Jets of distinguished curves}\label{2}

\subsection{The bundles of $\Cal C_A$--velocities}\label{2.1}
Let us recall the natural bundles $T^r_k$ of $r$th order
$k$--dimensional velocities on all smooth manifolds. By definition,
$T^r_kM=J^r_0(\Bbb R^k,M)$, so this is the bundle of $r$--jets of
parametrized $k$--dimensional (singular) submanifolds in $M$.  In
particular, $r$--jets of curves are elements in $T^r_1M$. The action
of all diffeomorphisms of $M$ on $T^r_kM$ is defined by jet
composition.  Let us consider a category of Cartan geometries of fixed
type $(G,P)$ and a class of generalized geodesics $\Cal C_A$, for a
$G_0$--invariant subset $A$ of $\frak n$. Then the jets of
distinguished curves of type $\Cal C_A$ form a natural subbundle
$T^r_{\Cal C_A}\subset T^r_1$ on parabolic geometries of type $(G,P)$.
Clearly, $T^r_{\Cal C_A}$ is a well defined functor, cf. Proposition
\ref{1.3}(2) above, however their values are not smooth bundles in
general, see the examples below. In the cases with $G_0$--invariant
subsets $A\subset \frak n$ we call the latter functors the
\textit{bundle of $r$th order velocities} of geodesics of type $\Cal
C_A$.

Our next goal is to prove that there always is a finite order $r$ for
which the entire geodesic is completely determined by a single value
in $T^r_{\Cal C_A}$.

\subsection{Jets of curves on $G/P$}\label{2.2}
Using Cartan's space $\S$, the development of curves defines a
bijection between smooth curves $c:I\to M$ defined on some
neighborhood $I$ of $0\in \Bbb R$ such that $c(0)=x_0$, and smooth
curves to $G/P$ which map $0$ to $o=eP$. Of course, this bijection is
compatible with taking jets in $x_0$, i.e.~two curves have the same
$\ell$--jet in $x_0$ if and only if the corresponding curves in $G/P$
have the same $\ell$--jet in $o$. By definition, this bijection also
respects generalized geodesics of any type. Thus to prove that
geodesics of some type $\Cal C_A$ are determined by some jet in one
point, it suffices to consider the homogeneous model $G/P$ and the
point $o$. We start by considering $A=\frak n$ (which of course
provides an estimate for any $A\subseteq\frak n$). Thus, we have to
study the curves $c^{b,X}(t) =b\exp(tX)P$, with $b\in P$ and $X\in
\frak n$, cf.  \ref{1.2}.

Since $b\exp(tX)=\exp(t\Ad_bX)b$ we see that
$c^{b,X}(t)=\exp(t\Ad_b\cdot X)P$. For any two curves $c(t)$ and
$d(t)$ in $G$, there is a uniquely determined curve $u(t)$ in $G$ such
that $c(t)=d(t)\cdot u(t)$. The projections of $c(t)$ and $d(t)$ to
$G/P$ coincide if and only if $u(t)\in P$ for all $t$. Thus the curves
$c^{b_1,X_1}$ and $c^{b_2,X_2}$ coincide if and only if the uniquely
determined curve $u$ such that
\begin{equation}
\exp(t\Ad_{b_1}X_1) =\exp(t\Ad_{b_2}X_2) \cdot u(t)
\end{equation}
has values in $P$. Since $\exp$ is analytic, the curve $u$ must be
analytic, too, and hence it has values in $P$ if and only if all
derivatives $u\d{i}(0) = (\tfrac d{dt})^iu(0)$ are tangent $P$. To
formulate this precisely, we use left logarithmic derivative $\delta
u:\Bbb R\to\fg$ of the curve $u:\Bbb R\to G$, see e.g.~\cite[p.
39]{KMS}.  In fact $\de u: T\Bbb R=\Bbb R\x\Bbb R\to \frak g$, $\de
u(t) = T\la_{u(t)^{-1}}\o T_tu$, but we shall identify the linear map
$\de u(t,\ ):\Bbb R\to \frak g$ with its value at the unit $1\in
T_t\Bbb R$.  Since knowing $\de u$ is equivalent to knowing $T u$, the
following Lemma is a simple observation.

\begin{lem*}
For each order $k\in\Bbb N$ we have $j^k_0\c1 =j^k_0\c2$ if and only
if the derivatives $(\de u)\d{i}(0)$ lie in $\p$ for all $i\leq k-1$.
\end{lem*}

\subsection{Some technicalities}\label{2.3}
In order to compute the derivatives of $\de u$ from formula
\ref{2.2}(1), we can use the Leibniz rule for the left logarithmic
derivative,
$$
\de(f\cdot g)(x) = \de g(x) + \Ad_{g(x)^{-1}} \de f(x),
$$
cf. \cite[p. 39]{KMS}, so it remains to compute the left
logarithmic derivative of the curve $t\mapsto \exp tX$. For later use,
we shall compute this expression with an arbitrary curve $Y:\Bbb R\to
\frak g$ instead of the line $tX$. By definition, the logarithmic
derivative $\de(f\o g)$ of the composition of two smooth maps $f:M\to
G$, $g:N\to M$ is given by $\de(f\o g)=(\de f)\o Tg$.  Thus, the key
ingredient is the formula for $\de(\exp): T\frak g\to \frak g$. The
proof of this formula for the right logarithmic derivative in \cite[p.
39]{KMS} can be easily adapted to our case, leading to 
$$
\de(\exp)(Y) = \sum_{p=0}^\infty \frac1{(p+1)!}\ad (-Y)^p.
$$ 
This proves:

\begin{lem*}
Let $Y:\Bbb R\to \frak g$ be a smooth curve with derivative $Y':\Bbb
R\to \frak g$. Then
$$
\de(\exp\o Y)(t)=\sum_{p=0}^\infty \frac1{(p+1)!}\ad(-Y(t))^p\cdot Y'(t)
.
$$ 
\end{lem*}

The first terms in the formula for $\de(\exp Y(t))$ read as
$$
Y'(t)-\tfrac12[Y(t),Y'(t)] + \tfrac16[Y(t),[Y(t),Y'(t)]] +\dots .
$$
Notice that if $Y$ has values in $\frak n$, then also $Y'$ has
values in $\frak n$, and compatibility of the grading of $\frak g$
with the Lie bracket implies that at most $k$ of these terms may be
non--zero for $|k|$--graded $\frak g$. Thus, for example,
\begin{align*}
\de(\exp Y(t)) &= Y'(t), \mbox{\ if $k=1$},
\\
\de(\exp Y(t)) &= Y'(t)-\tfrac12[Y(t),Y'(t)],  \mbox{\ if $k=2$},
\\ 
\de(\exp Y(t)) &= Y'(t)-\tfrac12[Y(t),Y'(t)] + 
\tfrac16[Y(t),[Y(t),Y'(t)]],  \mbox{\ if $k=3$}. 
\end{align*}
On the other hand, if $Y(t)=\ph(t)Y$ for some fixed $Y\in\fg$ and a smooth
function $\ph$, then $[Y(t),Y'(t)]=0$ and hence we always get
\begin{equation}
\de (\exp \ph(t)Y) = \ph'(t)Y.
\end{equation}

Applying the left logarithmic derivative to equation \ref{2.2}(1) yields
\begin{equation}
\de u(t) = \Ad_{b_1}X_1 - \Ad_{u(t)^{-1}}\Ad_{b_2}X_2.
\end{equation}
In particular, $\de u(0) = \Ad_{b_1}X_1 - \Ad_{b_2}X_2$, and this lies
in $\p$ if and only if $\Ad_{b_1}X_1$ and $\Ad_{b_2}X_2$ represent the
same class in $\frak g/\frak p$, i.e.~if the curves have the same
tangent vector at $0$.

Differentiating equation (2) at zero we obtain 
$$
(\de u)'(0) = -\ad_{(-u'(0))}\Ad_{b_2}X_2 = [u'(0),\Ad_{b_2}X_2],
$$ 
and $u'(0)$ is the image of $1\in T_0\Bbb R$ by $\de u(0)$.
Substituting (2) yields $(\de u)'(0) = [\de u(0),\Ad_{b_1}X_1]$.
Surprisingly, there is a general formula for $(\de u)^{(i)}(t)$ for
all $t\in\Bbb R$ and all orders $i$:

\begin{lem}\label{2.4} For all $i\ge1$,
$(\de u)^{(i)}(t) = (\ad(-\Ad_{b_1}X_1))^{i}(\de u(t))$.
\end{lem}

\begin{proof}
  Let us start with the first order derivative, so we have to prove
  $(\de u)'(t)= [\de u(t),\Ad_{b_1}X_1]$.  To do this, we have to
  compute the derivative of $t\mapsto \Ad_{u(t)^{-1}}:\R\>GL(\g)$.
  Clearly, $\tfrac d{dt}(t\mapsto\Ad_{u(t)^{-1}}) =(T\Ad \o
  T\nu)(u'(t))$, where $\nu$ is the inversion in $G$ and $T_t
  u=u'(t)$.  First, we will express $T_g\nu$ and $T_g\Ad$ in general.

{From} $\rh_g\o\nu\o\la_g=\nu$ we have 
$T_{g^{-1}}\rh_g \o T_g\nu \o T_e\la_g = T_e\nu$, thus 
$T_g\nu = -T_e\rh_{g^{-1}} \o T_g\la_{g^{-1}}$.
Similarly, $\Ad\o\la_g =\Ad_g\o\Ad$ implies 
$T_g\Ad\o T_e\la_g =\Ad_g\o T_e\Ad$, so 
$T_g\Ad =\Ad_g\o\ad\o T_g\la_{g^{-1}}$.
Altogether,
$$
\tfrac d{dt}\Ad_{u(t)^{-1}} =(\Ad_{u(t)^{-1}}\o\ad\o T\la_{u(t)}) 
\o(-T\rh_{u(t)^{-1}}\o T\la_{u(t)^{-1}})(u'(t)).
$$
Since $\Ad_g=T_e(\la_g\o\rh_{g^{-1}})$ and 
$\de u(t)=T\la_{u(t)^{-1}}\o u'(t)$ the latter expression 
equals  $(-\Ad_{u(t)^{-1}}\o\ad\o\Ad_{u(t)})(\de u(t))$. 
Thus,
$$
(\de u)'(t) =\Ad_{u(t)^{-1}}[\Ad_{u(t)}\de u(t),\Ad_{b_2}X_2] 
            =[\de u(t),\Ad_{u(t)^{-1}}\Ad_{b_2}X_2]
$$
and substituting $\Ad_{u(t)^{-1}}\Ad_{b_2}X_2=\Ad_{b_1}X_1- \de u(t)$
from \ref{2.3}(2) the claim follows. 

Now, let $i>1$ and assume that the formula is valid for all orders
less then $i$. Then
$$
(\de u)^{(i)}(t)=\ddt t(\ad(-\Ad_{b_1}X)^{(i-1)}\de u(t))
$$
and since $\ad(-\Ad_{b_1}X)^{(i-1)}$ is a linear map and we have
computed $(\de u(t))'$ already, we arrive at
$$
(\de u)^{(i)}(t) = \ad(-\Ad_{b_1}X)^{(i-1)}(\de u(t))'= 
\ad(-\Ad_{b_1}X)^{(i)}\de u(t),
$$
which is the required formula.
\end{proof}

Let us notice that we have also derived
the more general formula for the derivative of
$\Ad_{u(t)^{-1}}Y(t)$ with $Y: \Bbb R\to \frak n$. From the proof
above we conclude
\begin{equation}
\ddt t(\Ad_{u(t)^{-1}}Y(t)) = \Ad_{u(t)^{-1}}Y'(t) - [\de u(t),
\Ad_{u(t)^{-1}}Y(t)]
.\end{equation} 

As a simple consequence of this Lemma, we can prove that any geodesic
is determined by a finite jet in one point:

\begin{prop}\label{2.5}
Let $\g$ be a $|k|$--graded Lie algebra, and let $A\subseteq\frak n$ be
any $G_0$--invariant subset. If two geodesics of type $\Cal C_A$ have
the same $(k+2)$--jet in one point, then they coincide. 
\end{prop}
\begin{proof}
  As we have noticed in \ref{2.2} it suffices to consider $A=\frak n$,
  an we can complete the proof by showing that two curves $\c1$ and
  $\c2$ coincide if they have the same $(k+2)$--jet in $0$.  Denoting
  by $u:\Bbb R\to G$ the curve determined by equation \ref{2.2}(1),
  Lemma \ref{2.4} tells us that $(\de u)^{(i)}(0) =
  (\ad(-\Ad_{b_1}X_1))^{i}(\de u(0))$. By Lemma \ref{2.2}, the
  assumption on the $(k+2)$--jet in $0$ implies that
  $\ad(-\Ad_{b_1}X_1)^{i}(\de u(0))\in\frak p$ for all $i\leq k+1$.
  Since $b_1\in P$, we may hit this element with $\Ad_{b_1}^{-1}$, and
  the result remains in $\frak p$. Putting $X=X_1\in\frak n$ and
  $Z=\Ad_{b_1^{-1}}\de u(0) \in\p$ we conclude that
  $\ad(-X)^i(Z)\in\frak p$ for all $i=1,\dots,k+1$. Since $Z\in\frak
  p=\fg_0\oplus\dots\oplus\fg_k$ and $-X\in\frak
  n=\fg_{-k}\oplus\dots\oplus\fg_{-1}$, compatibility of the bracket
  with the grading implies that $\ad(-X)^i(Z)\in\frak
  g_{-k}\oplus\dots\oplus\frak g_{k-i}$. Putting $i=k+1$, we see that
  $\ad(-X)^{k+1}(Z)$ has to lie both in $\frak n$ and in $\frak p$, so
  it must be zero. This implies that $\de u^{\ell}(0)=0\in\frak p$ for
  all $\ell>k+1$, and thus $\c1=\c2$ and the claim follows.
\end{proof}
Let us remark at this point that the estimate $r=k+2$ on the jet
needed to pin down a geodesic is not at all sharp and we will improve
it heavily depending on a particular choice of the class of geodesics.

\subsection{Distinguished curves in a given direction}\label{2.5a}
The most natural way to approach the problem of distinguished curves
usually is to fix a point $x\in M$ and a tangent vector $\xi\in T_xM$,
and look for geodesics emanating from $x$ in direction $\xi$. Given a
$G_0$--invariant subset $A\in\frak n$, the basic question then is how
many geodesics of type $\Cal C_A$ pass through $x$ in direction $\xi$.
Of course, it may happen that there are no such geodesics. As before,
one may restrict the discussion to the point $o$ in the homogeneous
model $G/P$. Since the above question is perfectly geometric, the
answer for a tangent vector $\xi\in T_o(G/P)\cong\frak g/\frak p$ will
only depend on the $P$--orbit of $\xi$. Clearly, there is at least one
geodesic of type $\Cal C_A$ in direction $X$, if the image of $A$ in
$\frak g/\frak p$ meets the $P$--orbit of $\xi$. Otherwise put, if
$X\in\frak n\subset\fg$ is the unique element such that
$\xi=X+\frak p$, then there is at least one geodesic of type $\Cal
C_A$ in direction $\xi$ if $\Adb_b(X)\in A$ for some $b\in P$. 

Second, suppose that $A,B\subset\frak n$ are $G_0$--invariant subsets,
and that for each $X\in A$ there is an element $b\in P$ such that
$\Ad_bX\in B$, and vice versa. (Of course, this is a very restrictive
condition, since we are using $\Ad_b$, which does not leave $\frak n$
invariant, but it happens in interesting cases.) Then this gives rise
to a bijection between the sets $\Cal C_A$ and $\Cal C_B$ of curves in
$G/P$, and consequently, geodesics of type $\Cal C_A$ coincide with
geodesics of type $\Cal C_B$. 

Fix a $G_0$--invariant subset $A\subseteq\frak n$ and an element $X\in
A$, and consider the tangent vector $\xi=X+\frak p\in T_o(G/P)$.
Clearly, $c^{e,X}(t)=\exp(tX)$ is a geodesic of type $\Cal C_A$ in
direction $X$, and any other geodesic of that type can be written as
$c^{b,Y}$ with $b\in P$ and $Y\in A$. It is a general fact, see
\cite[2.10]{CS} that there are unique elements $b_0\in G_0$ and
$Z\in\frak p_+$ such that $b=b_0\exp(Z)=\exp(\Ad_{b_0}Z)b_0$. From the
definition of distinguished curves, we conclude that 
$$
c^{b_0\exp Z,Y} =c^{\exp(\Ad_{b_0}Z),\Ad_{b_0}Y},$$
and $\Ad_{b_0}Y\in A$. Hence any geodesic of type $\Cal C_A$ may be
written as $c^{\exp(Z),Y}$ for $Z\in\frak p_+$ and $Y\in A$. Hence we
conclude that the set of geodesics of type $\Cal C_A$ in direction
$\xi=X+\frak p$ can be equivalently described as 
$$
\{c^{\exp(Z),Y}:Z\in\frak p_+,Y\in A,\Adb_{\exp(Z)}\cdot Y=X\}. 
$$

Passing to a general curved geometry via developments as before, we
obtain 
\begin{prop*}
Let $(p:\Cal G\to M,\om)$ be a Cartan geometry of type $(G,P)$,
$x\in M$ a point, $\xi\in T_xM$ a tangent vector, and let
$A\subseteq\frak n$ be a $G_0$--invariant subset. Then there is a
geodesic of type $\Cal C_A$ through $x$ in direction $\xi$ if and only
if there are elements $u\in p^{-1}(x)\subset\Cal G$ and $X\in A$ such
that $\xi=T_up\cdot\om^{-1}(X)$. Moreover, for any such pair $(u,X)$,
one obtains a bijection between the set of geodesics of type $\Cal
C_A$ through $x$ in direction $\xi$ and the set
$\{c^{\exp(Z),Y}:Z\in\frak p_+,Y\in A,\Adb_{\exp(Z)}\cdot Y=X\}$ of
curves in $G/P$. This bijection is compatible with finite jets in $0$
in the obvious sense. 
\end{prop*}

Finally note that the curves $c^{\exp(Z_1),Y_1}$ and
$c^{\exp(Z_2),Y_2}$ have the same $\ell$--jet in $0$ respectively
coincide if and only if the same is true for $c^{e,Y_1}$ and
$c^{\exp(Z_1)^{-1}\exp(Z_2),Y_2}$, and we can write
$\exp(Z_1)^{-1}\exp(Z_2)$ as $\exp(Z)$ for some $Z\in\frak p_+$. Hence
we conclude that if for some $\ell$ and each $X\in A$ we can show that
any curve $c^{\exp(Z),Y}$ with $Y\in A$ which has the same $\ell$--jet
in $0$ as $c^{e,X}$ must actually equal $c^{e,X}$, then this implies
that any geodesic of type $\Cal C_A$ is uniquely determined by its
$\ell$--jet in a single point. 

\subsection{The $|1|$--graded case}\label{2.6}
For irreducible parabolic geometries we easily reach a complete
description. So we assume $\frak g=\frak g_{-1}\oplus\frak g_0\oplus
\frak g_1$ and $A=\frak n$. The main simplification in the
$|1|$--graded case comes from the fact that in this case $\frak p_+$
acts trivially on $\fg/\frak p$, so the $P$ action on this quotient
factorizes over $G_0$. In particular, for $Z\in\frak p_+=\fg_1$ and
$Y\in\frak n=\fg_{-1}$ we get $\Adb_{\exp(Z)}Y=Y$, so in view of
Proposition \ref{2.5a} it remains to compare the curves $c^{e,X}$ and
$c^{\exp(Z),X}$ with $Z\in\fg_1$. For the corresponding curve $u$, we
obviously get $\de u(0)=-[Z,X]-\tfrac12[Z,[Z,X]]$. For the two curves
having the same two--jet in $0$, we must have
$$
(\de u)'(0) = -[X_1,\de u(0)] =
[X_1,[Z,X_1]]+\tfrac12[X_1,[Z,[Z,X_1]]]\in\frak p,
$$
and thus $[X_1,[Z,X_1]]=0$. But this implies
$[X_1,[Z,[Z,X_1]]]=[Z,[X_1,[Z,X_1]]]=0$, and so $(\de u)^{(i)}(0)=0$
for all $i\ge 2$. Thus, we have proved:

\begin{prop*}
  Each generalized geodesic in an irreducible parabolic geometry is
  uniquely determined by its 2--jet in one point.
\end{prop*}

\subsection{The distinguished jets}\label{2.7}
Using the procedures from above, one may compute explicitly the jets
of all geodesics of type $\Cal C_A$. For the sake of simplicity, we
shall restrict ourselves again to the case of $|1|$--graded Lie
algebras. Thus, the value in $T^2_1(G/P)$ over the origin will always
determine a geodesic completely, and we shall compute explicitly the
algebraic description of the standard fibers of $T^2_{\Cal C_A}$.
Understanding the higher jets of geodesics is an interesting problem,
however the computations grow quickly out of hand.

Let us describe all distinguished curves in normal coordinates through
the origin, i.e.~we have to represent each geodesic in the form
$t\mapsto \exp(Y(t))P$ for a smooth curve $Y:\Bbb R\to\frak g_{-1}$
with $Y(0)=0$. This means, that rather than with formula \ref{2.2}(1),
we have to deal with
$$
\exp(Y(t))\cdot u(t) = \exp (t\Ad_{\exp Z}X)
$$
for $Z\in \frak g_1$ and $X\in A\subseteq \frak g_{-1}$.

Using the results in \ref{2.3} and formula \ref{2.4}(1),
straightforward computations yield
\begin{align*}
\de u(t) &= X + [Z,X] +\tfrac12[Z,[Z,X]]-\Ad_{u(t)^{-1}}(Y'(t))
\\
(\de u)'(t) &= \left[X+[Z,X] +\tfrac12[Z,[Z,X]],\Ad_{u(t)^{-1}}Y'(t)\right] -
\Ad_{u(t)^{-1}}Y''(t)
.\end{align*} 
The requirement $(\de u)^{(i)}(0)\in \frak p$, for  $i=0,1$
immediately implies 
\begin{align*}
Y'(0) &= X
\\
Y''(0) &= [X,[X,Z]]
.\end{align*}

Now it is easy to describe the standard fiber of $T^2_{\Cal C_A}$ as
follows.  The standard fiber of $T^2_1$ is the smooth manifold
$J^2_0(\R,\fg_{-1})_0$, which is naturally identified with
$\g_{-1}\x\g_{-1}$.  Hence the standard fiber of $T^2_{\Cal C_A}$ is a
subset in $\g_{-1}\x\g_{-1}$, which we have computed to be
$$
S=\left\{\begin{pmatrix}X\\ [X,[X,Z]]\end{pmatrix}: X\in A,Z\in\g_1\right\}.
$$
Recall that $A$ is assumed to be $G_0$--invariant, but not
necessarily a linear subspace. A good examples in which it is not a
subspace is given by the null cone in $\Bbb R^{p+q}$ in the setting of
Example \ref{1.5}(3). In that case, $[X,[Z,X]]$ happens to be a
multiple of $X$ for each $Z$, which corresponds to the fact that
geodesics in null directions are conformally invariant up to
parametrization. 
 
For every parabolic geometry of type $(G,P)$, there is the standard
embedding $i:P\>G^2_m=\operatorname{inv}J^2_0(\Bbb R^m,\Bbb R^m)_0$,
see e.g. \cite{Och, Sha}.  Further, the action of the structure group
$G^2_m$ on $J^2_0(\R,\R^m)_0$ transforms to the action on
$\g_{-1}\x\g_{-1}$, whose restriction to the subgroup $i(P)$ keeps the
subset $S$ invariant because the set $\Cal C_A$ of all geodesics is
$P$--invariant.

In fact, the action of $G_0$ obviously is the product of the adjoint
actions on $\frak g_{-1}\x \frak g_{-1}$, while the action of
$P_+=\exp\frak g_1$ comes by the very definition of the curves from
the left shift by the elements $\exp W$, $W\in \frak g_1$. Since
$\frak g_1$ is an abelian subalgebra, the action by $\exp W$ is given
by 
$$
\exp W\cdot \begin{pmatrix}Y'\\Y''\end{pmatrix} = 
\begin{pmatrix}Y'\\Y'' + [Y',[Y',W]]\end{pmatrix}.
$$
Hence we obtain an alternative description of the standard fiber as
the $P$--orbit of the $G_0$--invariant subspace $A\x\{0\}\subseteq
\frak g_{-1}\x\fg_{-1}$.

\section{Reparametrizations}\label{3}

In this section we shall generalize our basic question to: {\em When
  are two distinguished curves equal up to a change of
  parametrization?} Thus we shall discuss the non--parametrized
geodesics together with their preferred parametrizations.

\subsection{Technicalities}\label{3.1}
In order to deal with this question, we have to modify our basic equation
\ref{2.2}(1). The answer is positive if and only if there exist 
mappings $u:\R\>P$ and $\ph: \Bbb R\to \Bbb R$ such that 
\begin{equation}
\exp(\ph(t)\Ad_{b_1}X_1) =\exp(t\Ad_{b_2}X_2)\cdot u(t),
\end{equation}
where $\ph$ is a local reparametrization, i.e. we require
$\ph'(t)\neq 0$ and, for simplicity, $\ph(0)=0$. As discussed in the Remark
\ref{2.5a}, we may restrict ourselves to $G_0$--invariant subsets $A$,
$b_1=e$, $b_2=\exp Z$ with $Z\in \frak p_+$. 

The left logarithmic derivative of (1) gives, cf. \ref{2.3}(1)
\begin{equation}
\de u(t) = \ph'(t)X_1 - \Ad_{u(t)^{-1}}\cdot\Ad_{\exp Z}X_2.    
\end{equation}
In particular, $\de u(0)\in\p$ if and only if tangent vectors of the two 
distinguished curves at 0 are equal up to a scalar multiple.

By  formula \ref{2.4}(1), and the above equation (2), we get 
\begin{equation}
(\de u)'(t) = \ph''(t)X_1 -\ph'(t)[X_1,\de u(t)].
\end{equation}

Now, similarly as in the parametrized case we prove a general iterative
formula for $(\de u)^{(i)}$:

\begin{lem}\label{3.2}
For all $i\geq 1$ and at every $t\in \Bbb R$, with the notation as above
$$
(\de u)^{(i)} = \ph^{(i+1)}X_1 +
\sum_{k=1}^i(-1)^k\bigl(\sum_{{\bf j},{\bf a}}c_{{\bf j}, {\bf a}}
(\ph^{(j_1)})^{a_1}\!\!\dots(\ph^{(j_s)})^{a_s}\bigr)
(\ad_{X_1})^k(\de u)
$$
where the internal sum runs over all $s$--tuples of 
natural numbers ${\bf j}=(j_1,\dots,j_s)$, $j_1< j_2<\dots<j_s$, 
and $s$--tuples of arbitrary
natural numbers ${\bf a}=(a_1,\dots,a_s)$ such that
$a_1j_1+\dots + a_sj_s=i$ and $a_1+\dots+a_s=k$, and the coefficients 
$c_{{\bf j}, {\bf a}}$ are
$$
c_{{\bf j}, {\bf a}} = \frac{i!}{(j_1!)^{a_1}\dots (j_s!)^{a_s}a_1!\dots a_s!}
.$$
\end{lem}
\begin{proof}
In the case $i=1$, the entire sum in the formula has
just one possible term for $k=1$, $j_1=1$ and $a_1=1$. 
As we have seen, this is the correct formula (3).
The general case is proved by a tedious induction.
\end{proof}

\begin{remark*} As a hint for the induction mentioned in the proof above, 
let us describe in words, what the individual terms in the general 
formula mean. The value of $k$ says how many times $\ph$ occurs
in the term in question (and so many times $X$ hits $\de u$ via the adjoint
action and the sign is set appropriately), 
while the coefficients $c_{\mathbf j, \mathbf a}$ 
express in how many different
ways we may split $i$ derivatives onto $k$ copies of 
$\ph$'s in order to achieve the result 
$(\ph^{(j_1)})^{a_1}\!\!\dots(\ph^{(j_s)})^{a_s}$. 
Now, the differentiation of this formula and substitution from \ref{3.1}(3) 
means that we perform the last
derivative on one of the $\ph$'s in the individual terms in the formula, 
or we attach a new $\ph$  to the existing
terms which is differentiated only once.
But this is exactly how all splittings of $i+1$ (distinguishable) hits of $k$
(indistinguishable) targets are obtained from the answers to the same
question for $i$ derivatives and $k$ or $k-1$ targets. 
Either the last hit has been to some existing one among $k$ targets, 
i.e. we use the answer with $i$ hits and $k$ targets, 
or we have had to introduce a new target which was hit
once, i.e. we used the answer with $i$ hits and $k-1$ targets.
\end{remark*}

It is probably hard to deduce general results for all parabolic geometries
and all classes of distinguished curves from this formula, but let us see
how to use it in more specific situations.
 
\subsection{Irreducible parabolic geometries}\label{3.3}
We are going to give a complete answer 
to our question for $|1|$--graded algebras $\frak g$.
In order to decide when two distinguished paths $\c1$, $\c2$ 
parametrize the same curve we have to compute explicitly the 
consequences of $(\de u)\d{i}(0)\in\p$ 
in relation to the necessary and sufficient conditions for the solution of 
the given problem.
At the same time we shall get a complete and explicit description of the 
reparametrizations.

\begin{lem*} With the notation as above, $\de u(0)\in\p$
if and only if 
\begin{equation}
\ph'(0)X_1 =X_2.
\end{equation}
If $\de u(0)\in\p$, then $(\de u)'(0)\in\p$ if and only if
\begin{equation} 
\frac{\ph''(0)}{\ph'(0)^2}X_1 =  [X_1,[X_1,Z]],
\end{equation}
and if $i\ge2$ and $(\de u)^{(j)}(0)\in\p$ for all $j<i$, then
$(\de u)^{(i)}(0)\in\p$ if and only if
\begin{equation} 
\ph^{(i+1)}(0)
=\frac{(i+1)!}{2^i}\frac{\ph''(0)^i}{\ph'(0)^{i-1}},\mbox{\ for all $i\ge
2$}.
\end{equation}
\end{lem*}
\begin{proof} 
Since our algebra $\frak g$ is $|1|$--graded, all iterated adjoint actions
by $X_1$ on $\de u(0)$ vanish if the order is more then two. Thus only terms
with $k\le 2$ in Lemma \ref{3.2} may survive and the general 
formula for $i\ge 1$ reads
\begin{align*}
(\de u)^{(i)}(0) = {}&\ph^{(i+1)}(0)X_1 - \ph^{(i)}(0)[X_1,\de u(0)] +
{}\\
&\frac12
\sum_{\ell=1}^{i-1}\tfrac {i!}{\ell!(i-\ell)!}\ph^{(\ell)}(0)\ph^{(i-\ell)}(0)
\bigr)[X_1,[X_1,\de u(0)]].
\end{align*}
Indeed, this can be either proved by inserting into 
the general formula from Lemma \ref{3.2} or directly by induction.  

Next, recall 
$\de u(0)= \ph'(0)X_1- X_2-[Z,X_2]-\frac12[Z,[Z,X_2]]$. Thus (1) is
obvious. We shall assume $\de u(0)\in \frak p$ and therefore
$$
\de u(0)= -\ph'(0)([Z,X_1]+\tfrac12[Z,[Z,X_1]]).
$$
Now, (2) follows from the general formula with $i=1$. The most interesting
step is the case $i=2$ (i.e. we deal with the third order jets of the curves,
so that these must be determined by the lower order derivatives already).
Indeed, the substitution of equalities in (1) and (2) into 
the general formula yields
\begin{align*}
(\de u)^{(i)}(0) = {}&\ph^{(i+1)}(0)X_1+\ph^{(i)}\ph'(0)[X_1,[Z,X_1]] -
\\
&{}\frac14\sum_{\ell=1}^{i-1}\begin{pmatrix}i\\\ell\end{pmatrix}
\ph^{(\ell)}(0)\ph^{(i-\ell)}(0)\ph'(0)[X_1,[X_1,[Z,[Z,X_1]]]] + 
\mbox{term in $\frak g_0$}
\\
= {}&\biggl(\ph^{(i+1)}(0) -\frac{\ph^{(i)}(0)\ph''(0)}{\ph'(0)} -\frac14
\sum_{\ell=1}^{i-1}\begin{pmatrix}i\\\ell\end{pmatrix}
\ph^{(\ell)}(0)\ph^{(i-\ell)}(0)\frac{\ph''(0)^2}{\ph'(0)^3}\biggr)X_1  
\\
&{}+ \mbox{ term in $\frak g_0$.}
\end{align*}
The structure of the latter equation implies that $\ph^{(i+1)}(0)$ is
determined uniquely in terms of the values $\ph^{(k)}(0)$ with $k\le i$ and
a direct computation checks that the formula in (3) is correct.  
\end{proof}

Let us summarize what we have achieved so far.  If the conditions of
(1) and (2) are satisfied, than $\ph(0)$, $\ph'(0)$ and $\ph''(0)$ are
determined by the choice of the tangent vectors to the curve and by
the element $Z\in\frak g_1$ and we may define all other derivatives of
$\ph$ by the formula (3).  In particular, the special case $i=2$
yields
\begin{equation}
\ph'''(0)=\frac32\frac{\ph''(0)^2}{\ph'(0)}
\end{equation}
which reminds the well known Schwartzian differential equation. We shall see
that the formulae for $\ph^{(i)}(0)$ determine an analytic local solution
for this equation.
 
If we denote $a=\ph'(0)$ and $b=\ph''(0)$, the Taylor development of 
the function $\ph$ at $0$ must be 
$$
\ph(t) = at +\frac1{2}bt^2 +\frac1{6}\frac{3}{2}\frac{b^2}{a}t^3 +\cdots
                +\frac1{(i+1)!}\frac{(i+1)!}{2^i}\frac{b^i}{a^{i-1}}t^{i+1} 
+\cdots.
$$
Thus, we have obtained the geometric series 
$\ph(t) = at\sum_{i=0}^{\infty}(\tfrac{bt}{2a})^i$
which converges locally around $0$ and its value is 
$$
\ph(t) = at(1-\frac{b}{2a}t)^{-1}
.$$

If we want to allow the reparametrizations with $\ph(0)\neq0$, 
we have just to replace equation \ref{3.1}(1) by
$
\exp(\ph(t)-\ph(0))\Ad_{b_1}X_1 =\exp(t\Ad_{b_2}X_2)\cdot u(t)
$
and the result differs only by adding the value 
$\ph(0)$ to the fraction above.
In such a case the reparametrization takes a form 
$$
\ph(t) =\frac{At+B}{Ct+D},\mbox{\ where\ } 
A=\ph'(0)-\frac{\ph''(0)}{2\ph'(0)}\ph(0),\ 
B=\ph(0),\ C=-\frac{\ph''(0)}{2\ph'(0)},\ D=1
.$$ 
In particular, the solution with $\ph''(0)=0$ yields the affine
reparametrization of the curve which of course have to be geodesics,
too. The determinant of the matrix $\smatrix{A&B\\ C&D}$ is
$\ph'(0)\neq 0$, so we may normalize this to just 1 and we have proved:

\begin{prop}\label{3.4}
  Suppose that $\fg$ is $|1|$--graded.  If the curves $c^{b_1,X_1}$
  and $c^{b_2,X_2}$ coincide as unparametrized curves, then the
  corresponding local reparametrization $\ph$ has the form
  $\ph(t)=\frac{At+B}{Ct+D}$, where $\smatrix{A&B\\C&D} \in SL(2,\R)$.
  Conversely, if $c=c^{b,X}$ is a parametrized geodesic then all
  curves $c\o\ph$ with reparametrizations $\ph:\Bbb R\to \Bbb R$ of
  the latter form are again geodesics
if and only if there is $Z\in \frak g_1$ such that
$[X,[X,Z]]=X$.
\end{prop}
\begin{proof}
It remains to prove the second statement. Obviously we may restrict
ourselves to the case when $\ph(0)=0$. Then each $\ph$ satisfies
all conditions from Lemma \ref{3.3}, provided there is
a suitable $Z$ for (2).
\end{proof}

Reparametrizations of the above type are called projective, 
see \cite{BE-conf}, where they are obtained as solutions of the Schwartzian 
differential equation $\ph'''=\tfrac3{2}\tfrac{(\ph'')^2}{\ph'}$.

\begin{kor*} 
  Suppose that $\fg$ is $|1|$--graded. Then the curves $c^{e,X_1}$ and
  $c^{\exp Z,X_2}$ parametrize the same unparametrized geodesic if and
  only if there are $a\neq 0$ and $b$ such that $X_2=aX_1$ and
  $[X_2,[X_2,Z]]=bX_1$.  This is equivalent to the existence of the
  projective local reparametrization $\ph$ which is uniquely
  determined by the initial condition $\ph(0)=0$, $\ph'(0)=a$, and
  $\ph''(0)=b$.
\end{kor*}

\begin{example}
  In the following examples we use the obvious fact that in the case
  of a $|1|$--grading, elements of $P$ of the form $\exp(Z)$ for
  $Z\in\fg_1$ act trivially on $T_o(G/P)=\g/\p$. Then the $P$--action
  on this space factorizes over $G_0$.

\noindent
(1) Conformal Riemannian structures correspond to $G=O(p+1,q+1)$ and
the parabolic subgroup $P$ as in \ref{1.5}(3). In an appropriate
matrix representation, the grading of the Lie algebra $\g$ has the form 
\begin{gather*}
\g_{-1}=\left\{\smatrix{0&0&0\\X&0&0\\0&-X^tJ&0}:X\in\R^{p+q}\right\},\ 
\g_0=\left\{\smatrix{a&0&0\\0&A&0\\0&0&-a}:A\in\frak
o(p,q),a\in\R,\right\},
\\ 
\g_1=\left\{\smatrix{0&Z&0\\0&0&-JZ^t\\0&0&0}:Z\in\R^{p+q*}\right\}.
\end{gather*}
Here $J$ is the matrix defining the standard pseudo--metric of signature 
$(p,q)$ on $\R^{p+q}=\g_{-1}$.

A direct calculation shows that $[X,[X,Z]]=-2Z(X)X-||X||^2JZ^t$, where
$||X||^2=X^tJX$ and $Z(X)=ZX$ is a real number.  Obviously, the space
$\g_{-1}$ splits into three different orbits of the action of $G_0$
according to the sign of $\|X\|^2$.  The orbit of null--vectors is of
particular interest, since $[X,[X,Z]]=-2Z(X)X$ in that case.  This
just means that all distinguished curves with the common tangent
null--vector differ by a reparametrization which recovers the
classical result that the null geodesics of the metrics in the
conformal class together with the class of projective parametrizations
are invariants of the conformal structure. Of course, these curves
will have their tangent vectors null in all their points.

For all tangent vectors which are not null, the second derivative may
be chosen arbitrarily. So that the standard fiber $S$ in \ref{2.7} has
arbitrary entries in the bottom row if $X$ is not null, but only
multiples of $X$ if $X$ is null. On the other hand, there always is an
element $Z\in\fg_1$ such that $[[Z,X],X]=X$, so all geodesics carry a
natural projective structure.

(2) Almost Grassmannian structures.  In this case, $G=SL(n+m,\Bbb R)$
and the parabolic subgroup $P$ is the stabilizer of $\Bbb
R^n\subset\Bbb R^{n+m}$, so it consists of block upper triangular
matrices with two blocks of sizes $n$ and $m$. On the infinitesimal
level,
\begin{gather*}
\g_{-1}=\{\smatrix{0&0\\X&0}:X\in\R^{mn}\},\  
\g_0=\{\smatrix{A&0\\0&B}:\text{tr}(A)+\text{tr}(B)=0\},
\\ 
\g_1=\{\smatrix{0&Z\\0&0}:Z\in\R^{nm}\}.
\end{gather*}
First, it is easy to see that the subgroup $G_0$ consists of block
diagonal matrices, and its action on $\g_{-1}$ is given by $X\mapsto
TXS^{-1}$, $(S,T)\in G_0$. Thus two elements of $\fg_{-1}$ lie in the same
$G_0$--orbit if and only if they have the same rank. Further, the
computation of the iterated bracket yields $[X,[X, Z]]=-2XZX$. In
particular, the choice of the pseudoinverse matrix $Z=X^\dagger$
provides always a multiple of $X$ and so all generalized geodesics enjoy the
distinguished projective structure. If the rank of $X$ is one, then we
may choose $X$ to be the matrix with the left upper element $x_{11}=1$
and all other $0$.  Then $[X,[X,Z]]$ equals to $z_{11}X$ for all $Z$
and so this behavior must be shared by all matrices of rank one.
Thus, the directions corresponding to rank one matrices behave like
null directions in pseudo--conformal geometries. The other extreme is
that $X$ has maximal rank. Then one gets a lot of freedom in the
available second derivatives of the curves. The case that all elements
of $\fg_{-1}$ are possible second derivatives occurs only if $m=n$ and
$X$ has rank $n$. 

(3) Projective structures are the special case $n=1$ of Example (2)
above. In this case, the rank of $X\ne 0$ is always one. More
explicitly, the product $ZX$ is a real number, so the bracket
$[X,[X,Z]]$ is always a multiple of $X$.  {From} this it follows that
all unparametrized distinguished curves are determined by the
direction in a given point. This agrees with the classical definition
of a projective structure as a class of affine connections sharing the
same unparametrized geodesics.  All such connections are parametrized
by smooth one--forms on the base manifold and they correspond to the
Weyl connections defined in \cite{CS-weyl}.
\end{example}

\section{More refinements}\label{4}
In this section we improve the estimates on the jet in a point needed
to pin down a geodesic for geodesics of certain types. The most
general results is Theorem \ref{4.3} but since the proofs of these
results are a bit technical, we prefer to discuss two simpler special
cases first. 

\subsection{Curves tangent to $T^{-1}M$}\label{4.1} 
Let $M$ be any manifold equipped with a parabolic geometry of some
fixed type $(G,P)$. A (generalized) geodesics with development of the
form $c^{b,X}$ emanates in a direction in $T^{-1}M$ if and only if
$X\in\g_{-1}$. Thus we are dealing with distinguished curves of type
$\Cal C_{\frak g_{-1}}$ and from Proposition \ref{1.3} we see that
they will be tangent to the distribution $T^{-1}M$ in all points.

To discuss geodesics of type $\Cal C_{\fg_{-1}}$, by Proposition
\ref{2.5a} we have to fix $X\in\fg_{-1}$ and study the curves
$c^{\exp(Z),Y}$ for $Z\in\frak p_+=\fg_1\oplus\dots\oplus\fg_k$ and
$Y\in\fg_{-1}$ such that $\Adb(\exp(Z))(Y)=X$. Since $Y\in\fg_{-1}$ we
get $\Adb(\exp(Z))(Y)=Y$ for any $Z$, so we have to consider all
curves of the form $c^{\exp(Z),X}$ with $Z\in\frak p_+$. By
\cite[2.10]{CS} we get a nicer presentation of $\exp(Z)$. Namely, there
are unique elements $Z_i\in\fg_i$ for $i=1,\dots,k$ such that
$\exp(Z)=\exp(Z_1)\cdots\exp(Z_k)$. Since $\Ad(\exp(W))=e^{\ad(W)}$
for each $W\in\fg$ we get 
\begin{equation*}
\Ad_{\exp Z}X = \sum_{i_1,\dots,i_k} \tfrac1{i_1!\cdots i_k!}
                                (\ad Z_1)^{i_1}\cdots(\ad Z_k)^{i_k}X.
\end{equation*}
Moreover, since $X\in\fg_{-1}$ a summand in the right hand side lies
in $\fg_\ell$ if and only if $i_1+2i_2+\dots+ki_k=\ell+1$. 

We need another observation for the proof: Suppose that $Y\in\fg$ is
any element. The Jacobi identity reads as
$\ad_X\o\ad_Y=\ad_{[X,Y]}+\ad_Y\o\ad_X$. Inductively, this implies
that $\ad_X^n\o\ad_Y$ can be written as a linear combination of terms
of the form $\ad_{\ad_X^i(Y)}\o\ad_X^j$ with $0\leq i,j$ and $i+j=n$.
In particular, if $\ad_X^{\ell+1}(Y)=0$ for some $\ell\geq 0$, then
for each $n>\ell$ there is a linear map $\ph$ such that
$\ad_X^n\o\ad_Y=\ph\o\ad_X^{n-\ell}$. Of course, it is not difficult
to compute $\ph$ explicitly, but we will not need this explicit form. 

\begin{prop*}
  A parametrized generalized geodesic of type $\C_{\frak g_{-1}}$ in a
  parabolic geometry corresponding to a $|k|$--grading of $\fg$ is
  uniquely determined by its $(k+1)$--jet in a single point.
\end{prop*}
\begin{proof}
  Of course, we have proved this for $k=1$ in \ref{2.6}. In view of
  the above discussion and the last observation in \ref{2.5a} we have
  to show that for each fixed $X\in\fg_{-1}$ any curve of the form
  $c^{\exp(Z_1)\cdots\exp(Z_k),X}$ with $Z_i\in\fg_i$ which has the
  same $k+1$--jet in $0$ as $c^{e,X}$ actually equals $c^{e,X}$. 

Given $Z_1,\dots,Z_k$ define
$W:=\Ad(\exp(Z_1)\cdots\exp(Z_k))(X)-X\in\frak p$. From the above
discussion we see that 
\begin{equation}
  W=\sum_{i_1,\dots,i_k} \tfrac1{i_1!\cdots i_k!}  (\ad
Z_1)^{i_1}\cdots(\ad Z_k)^{i_k}X,
\end{equation}
where the sum is over all $(i_1,\dots,i_k)$ such that
$0<i_1+2i_2+\dots+ki_k\leq k+1$. Considering the curve $u(t)$
associated to $c^{e,X}$ and $c^{\exp(Z_1)\cdots\exp(Z_k),X}$ by
equation \ref{2.2}(1), we see from \ref{2.3} that $\de u(0)=-W$ and
Lemma \ref{2.4} implies that $(\de
u)^{(i)}(0)=(-1)^{i+1}\ad_X^{i}(W)$. Consequently by Lemma \ref{2.2}
proving the result boils down to showing that $\ad_X^i(W)\in\frak p$
for all $i\leq k$ implies $\ad^i_X(W)\in\frak p$ for all $i\in\Bbb N$.

For each $\ell=1,\dots,k$ define $W'_\ell$ to be the sum of those
terms in the expression (1) for $W$ for which all $i_j$ with $j>\ell$
are zero, and put $W''_\ell=W-W'_\ell$. In particular, we have
$W''_k=0$, i.e.~$W'_k=W$. 

\noindent
\textbf{Claim}: If $\ad_X^i(W)\in\frak p$ for all $i\leq\ell$, then
for each $j\leq\ell$ we have $\ad_X^{j+1}(Z_j)=0$, and for each
$n>\ell$ we get $\ad_X^n(W'_\ell)\in\frak p$. 

We prove this claim by induction on $\ell$. If $\ell=1$, we know that
$\ad_X(W)\in\frak p$. Looking at formula (1) for $W$ and taking into
account that $X\in\fg_{-1}$ we see that $\ad_X(W)\in\frak p$ implies
(and is actually equivalent to) $[X,[Z_1,X]]=0$ and thus to
$\ad_X^2(Z_1)=0$. Hence it remains to show that $\ad_X^n(W'_1)\in\frak
p$ for all $n>1$. By definition,
$W'_1=\sum_{i=1}^{k+1}\tfrac{1}{i!}\ad_{Z_1}^iX$. Thus
$\ad_X^n(W'_1)\in\frak p$ is equivalent to $\ad_X^n\o\ad_{Z_1}^i X=0$
for $i\leq n$. From above we know that $\ad_X^2(Z_1)=0$ implies that
$\ad_X^n\o \ad_{Z_1}=\ph\o\ad_X^{n-1}$, so inductively we conclude
that $\ad_X^n\o \ad_{Z_1}^i=\ps\o \ad_X^{n-i+1}\o\ad_Z$ for some
linear map $\ps$ and by assumption $n-i+1>0$. Hence applying this
element to $X$ we get $\ps\o \ad_X^{n-i+2}(Z)$ which vanishes since
$n-i+2\geq 2$. This completes the proof of the case $\ell=1$.

Assume inductively that $\ell>1$ and we have proved the result for
$\ell-1$. Given that $\ad_X^i(W)\in\frak p$ for all $i\leq\ell$, we by
induction conclude that $\ad_X^{j+1}(Z_j)=0$ for
$j=1,\dots,\ell-1$. Moreover, we know by induction that
$\ad_X^\ell(W)\in\frak p$ implies $\ad_X^\ell(W''_{\ell-1})\in\frak
p$. By definition of $W''_{\ell-1}$ the only term in $
\ad_X^\ell(W''_{\ell-1})$ which does not automatically lie in $\frak
p$ is $\ad_X^\ell([Z_\ell,X])$, so we conclude that
$\ad_X^{\ell+1}(Z_\ell)=0$. Hence it remains to show that
$\ad_X^n(W'_\ell)\in\frak p$ for all $n>\ell$. Since we know by
induction that $\ad_X^n(W'_{\ell-1})\in\frak p$, it suffices to
consider $\ad_X^n(W'_\ell-W'_{\ell-1})$. Now from the expression (1)
for $W$ we conclude that 
$$
W'_\ell-W'_{\ell-1}=\sum_{i_1,\dots,i_{\ell}} \tfrac1{i_1!\cdots i_{\ell}!}
                                (\ad Z_1)^{i_1}\cdots(\ad
                                Z_\ell)^{i_\ell}X,
$$
with the sum going over $i_\ell>0$ and $i_1+2i_2+\dots+\ell
i_\ell\leq k+1$. Obviously, $\ad_X^n(W'_\ell-W'_{\ell-1})\in\frak p$
is equivalent to vanishing of
$\ad_X^n\o\ad_{Z_1}^{i_1}\o\dots\o\ad_{Z_\ell}^{i_\ell}$ for all
multi--indices $(i_1,\dots,i_\ell)$ such that $i_1+2i_2+\dots+\ell
i_\ell\leq n$. Since $\ad_X^{j+1}(Z_j)=0$, we see from above that
$\ad_X^m\o\ad_{Z_j}=\ph\o\ad_X^{m-j}$ for $m>j$. Inductively we
conclude that for $m>ji_j$ we get
$\ad_X^m\o\ad_{Z_j}^{i_j}=\ps\o\ad_X^{m-ji_j}$ for some linear map
$\ps$. Thus we conclude that 
$$
\ad_X^n\o\ad_{Z_1}^{i_1}\o\dots\o\ad_{Z_\ell}^{i_\ell}=
\tilde\ps\o\ad_X^{n-i_1-2i_2-\dots-\ell(i_\ell-1)}\o\ad_{Z_\ell},
$$
and by assumption $n-i_1-2i_2-\dots-\ell(i_\ell-1)\geq\ell$. Thus
applying the right hand side to $X$, we obtain $\ad_X^r(Z_\ell)$,
and by construction $r\geq\ell+1$, so this vanishes. Hence the proof
of the claim is complete. 

But taking the claim in the case $\ell=k$, we see that $\ad^i_X(W)=0$
for all $i\leq k$ implies that $\ad_X^n(W'_k)\in\frak p$ for all
$n>k$. Since we have observed above that $W'_k=W$, this completes the
proof. 
\end{proof}

\subsection{The case $A=\frak g_{-k}$}\label{4.2}
The other extreme class of geodesics on a manifolds $M$ equipped with
a parabolic geometry of type $(G,P)$ with $|k|$--graded $\frak g$ is
provided by the generalized geodesics of type $\Cal C_{\fg_{-k}}$. Of
course, for a point $x\in M$ and a tangent vector $\xi\in T_xM$ one
must have $\xi\in T_xM\setminus T^{-k+1}_xM$ in order to have a
nontrivial geodesic of type $\Cal C_{\fg_{-k}}$ in direction $\xi$. On
the other hand, this condition is not sufficient for such a geodesic,
and the directions of these geodesics usually form a smaller cone in
each tangent space. 

An important special case is parabolic contact geometries, i.e.~those
geometries corresponding to $|2|$--gradings, such that $\fg_{-2}$ has
dimension one and the bracket $\fg_{-1}\x\fg_{-1}\to\fg_{-2}$ is
non--degenerate. These geometries always have an underlying contact
structure. In these cases geodesics of type $\Cal C_{\fg_{-2}}$ always
exist for all directions in $TM\setminus T^{-1}M$.  A very well known
instance of this type of generalized geodesics is provided by the
Chern--Moser chains on hypersurface type CR--structures. A slightly
more general example of this type was studied for 6--dimensional
CR--structures of codimension 2, in \cite{SchmSl}.

Let us recall that reparametrizations of the form $\ph(t)
=\frac{At+B}{Ct+D}$ with $A\ne 0$ and $AD-BC=1$ are called projective.

\begin{thm*}
  Each generalized geodesic of type $\Cal C_{\g_{-k}}$ in a parabolic
  geometry of type $(G,P)$ corresponding to a $|k|$--grading on $\frak
  g$ is uniquely determined by its $2$--jet in a single point.
  Moreover, if two of such curves coincide up to parametrization, then
  this reparametrization is projective. Conversely, given a
  generalized geodesic of type $\Cal C_{\frak g_{-k}}$ corresponding
  to $(u,X)\in \Cal G\x\frak g_{-k}$, every projective change of
  parametrization defines a geodesic of the same type if and only if
  there exists a $Z\in\frak g_k$ such that $[X,[X,Z]]=X$.
\end{thm*}
\begin{proof}
  From \ref{2.5a} and \ref{4.1} we know that for each $X\in\fg_{-k}$
  we have to compare $c^{e,X}$ to all curves of the form $c^{b,Y}$
  with $b=\exp(Z_1)\cdots\exp(Z_k)$ for $Z_i\in\fg_i$, $Y\in\fg_{-k}$
  and $\Adb(b)(Y)=X$. The last condition immediately implies that
  $Y=X$. Expanding $W=\Ad(b)(X)-X$ as in equation \ref{4.1}(1), we
  conclude that if this expression has trivial component in 
$\fg_{-k+1}$, then
  $[Z_1,X]=0$. Hence we may omit all terms in the expansion for which
  $i_1$ is the only nonzero index.  Vanishing of the component in
  $\fg_{-k+2}$ then implies $[Z_2,X]=0$, so we may omit terms in which
  only $i_1$ and $i_2$ are nonzero.  Inductively, we get
  $[Z_\ell,X]=0$ for all $\ell=1,\dots,k-1$. Hence we conclude that
  $\de u(0)=-[Z_k,X]-\frac12[[Z_k,[Z_k,X]]$.  Now $(\de u)'(0)\in\p$
  implies $[X,[Z_k,X]]=0$ and so $(\de u)'(0)=0$ exactly as in
  \ref{2.6}.
  
  Concerning reparametrizations, we may adapt the proofs of Lemma
  \ref{3.3} and Proposition \ref{3.4} along the same lines. Using the
  notation from there, the condition $\de u(0)\in\p$ implies
  $X_2=\ph'(0)X_1$ and moreover $[Z_\ell,X_2]=0$ for all $\ell\leq
  k-1$, inductively as above, and this is the only difference to the
  $|1|$--graded case.  Further, $(\de u)'(0)\in\p$ if and only if
  $\ph''(0)X_1=\ph'(0)^2 [X_2,[X_2,Z_2]]$ and we finish the proof
  exactly as in the $|1|$--graded case.
\end{proof}

More generally, let us consider generalized geodesics of type $\Cal
C_{\frak g_{-j}}$ with arbitrary $j$. Geodesics of this type are
always curves with tangents in $T^{-j}M$ and they emanate from a given
point in $M$ in certain directions in $T^{-j}M\setminus T^{-j+1}M$.

\begin{thm}\label{4.3}
  Each generalized geodesic of type $\Cal C_{\g_{-j}}$ in a parabolic
  geometry of type $(G,P)$ with a $|k|$--graded $\frak g$, $1\le j\le
  k$ is uniquely determined by its $r$--jet in a single point provided
  that $rj\geq k+1$.
\end{thm}
\begin{proof}
  This is a combination of the proofs of Theorem \ref{4.2} and of
  Proposition \ref{4.1} with minor generalizations, so we just outline
  the basic steps: For $X\in\fg_{-j}$ we have to compare $c^{e,X}$ to
  $c^{b,Y}$ for $b=\exp(Z_1)\cdots\exp(Z_k)$ with $Z_i\in\fg_i$ and
  $Y\in\fg_{-j}$ and $\Adb(b)(Y)=X$. This immediately implies $Y=X$,
  and we put $W=\Ad(b)(X)-X$ and expand this as in \ref{4.1}(1). The
  proof boils down to showing that $\ad_X^i(W)\in\frak p$ for $i\leq
  r$ implies the same result for all $i$. As in \ref{4.2}, $W\in\frak
  p$ implies that $[Z_\ell,X]=0$ for all $\ell<j$, so in the notation of
  the proof of Proposition \ref{4.1} we obtain $W'_{j-1}=0$.
  
  The analog of the claim in the proof of Proposition \ref{4.1} is
  that if $\ad_X^i(W)\in\frak p$ for all $i\leq\ell$, then for each
  $s\leq\ell$ and $m<(s+1)j$, we get $\ad_X^{s+1}(Z_m)=0$, and further
  $\ad_X^n(W'_{j\ell-1})=0$ for all $n>\ell$. This is proved by
  induction using the same arguments as in \ref{4.1}.

  For $\ell=r$, we obtain $jr\geq k+1$, and as in \ref{4.1}, $W'_k=W$,
  and we conclude that $\ad_X^i(W)\in\frak p$ for all $i\leq r$ implies
  the same property for all $i$ as required. 
\end{proof}

The following two examples expose the diversity of the possible
behavior of various classes of distinguished curves in specific
parabolic geometries.  All claims may be checked by direct
computations following the results above and their more detailed
version may be also found in \cite{Zadnik-dissertation}.

\begin{example}\label{4.4}
  Let us briefly illustrate the general results in the simplest cases
  of parabolic contact structures, so we are dealing with
  $|2|$--gradings such that $\frak g_{-2}$ is one--dimensional and the
  bracket $\fg_{-1}\x\fg_{-1}\to\fg_{-2}$ is non--degenerate. As we
  have mentioned in \ref{4.2}, we get in each direction outside the
  contact subbundle geodesics of type $\Cal C_{\g_{-2}}$ which
  generalize the Chern--Moser chains for CR--structure. From \ref{4.3}
  we know that they are determined by their two--jet in a point as
  parametrized curves, and it follows that they are 
uniquely determined by their direction in
  one point up to parametrization, by dimension reasons. 
Moreover, each such geodesic
  carries a natural projective structure of distinguished
  parametrizations.
  
  Apart of these types of generalized geodesics, there are several
  other possibilities for non--equivalent types of geodesics as we may
  observe already at the simplest example of $G$ being a real form of
  $SL(3,\Bbb C)$ and $P$ the Borel subgroup:

\noindent
(1) $G=SL(3,\Bbb R)$.  The corresponding geometries are the
Lagrangian contact structures on 3--dimensional manifolds,
i.e.~three dimensional contact structures endowed with a decomposition
of the contact subbundle into a direct sum of two line subbundles,
cf. \cite{Tak}.  Geometrically, there are four different classes of
tangent vectors. First, we have vectors tangent to one of the two
subbundles (two classes), then there are the remaining vectors in the
contact subbundle, and finally those outside of the contact
subbundle. 

The subgroup $P$ consist of all elements of $G$ which are upper
triangular, so on the Lie algebra level, we obtain $\frak n$ as the
subalgebra of strictly lower triangular matrices, with the two entries
directly below the main diagonal corresponding to $\fg_{-1}$ and the
entry in the lower left corner corresponding to $\fg_{-2}$. The action
of the subgroup $G_0$ rescales each entry of a matrix in $\frak n$ by
a nonzero factor, so the $G_0$--orbits in $\frak n$ are determined
simply by the nonzero entries of a matrix. 

First, there are two canonical invariant subspaces in $\g_{-1}$ which
correspond to the Lagrange subspaces of the contact distribution.
They are $A_1=\left\{\smatrix{0&0&0\\\ast&0&0\\0&0&0}\right\}$ and
$A_2=\left\{\smatrix{0&0&0\\0&0&0\\0&\ast&0}\right\}$, respectively,
where the star denotes a nonzero entry. Generalized geodesics of these
types exist exactly in directions tangent to one of the two line
subbundles, so the two classes are disjoint but have the same
properties.  In both cases they behave just like null--geodesics in
conformal geometry, i.e.\ each such curve is determined by its 2--jet in one
point and with a given tangent vector there is a 1--dimensional family
of parametrized generalized geodesics determined by elements of the
form $\smatrix{0&\ast&0\\0&0&0\\0&0&0}$ and
$\smatrix{0&0&0\\0&0&\ast\\0&0&0}\in\p_+$, respectively.  Moreover,
all curves from this family coincide up to a projective
reparametrization.

For $A=\g_{-1}$ we get directions in the contact distribution. From
\ref{4.1} we know that such curves are determined by their $3$--jet in
one point. There is a 3--dimensional family of parametrized
generalized geodesics (corresponding to all elements in $\p_+$)
sharing a given tangent vector, which is not tangent to one of the two
line subbundles.  Admissible reparametrizations are the projective
ones, so the dimension of the space of unparametrized generalized
geodesics with the common direction in $T^{-1}M$ but outside of the
Lagrange subspaces is two.

Now we discuss the curves emanating in directions which do not belong
to the contact distribution.  For $A=\g_{-2}$ we obtain the 
analog of CR--chains as described in Theorem \ref{4.3}. 

Besides these chains, there are another curves going in all directions
except those in the contact distribution; this class of curves
corresponds to the generic choice of
$A=\left\{\smatrix{0&0&0\\\ast&0&0\\\ast&\ast&0}\right\}$.  Curves of
this type are determined by a 2--jet and to any tangent vector there
is a 3--dimensional family of generalized geodesics.  This set is
parametrized by elements of $\p_+$.  In contrast to the previous
cases, there are no two curves with the common tangent vector, which
would be the same up to a reparametrization. So here only affine
reparametrizations are allowed.

The two $G_0$--orbits in $\frak n$, which have not yet been mentioned
are $\left\{\smatrix{0&0&0\\\ast&0&0\\\ast&0&0}\right\}$ and
$\left\{\smatrix{0&0&0\\0&0&0\\\ast&\ast&0}\right\}$. Any element of
either of these can by mapped to $\fg_{-2}$ by some $\Ad_b$ with $b\in
P$, and vice versa. Hence from \ref{2.5a} we know that these lead to
the same curves as $A=\g_{-2}$, and thus the discussion is complete.

(2) $G=SU(2,1)$.  The corresponding geometries are non--degenerated
strictly pseudoconvex $3$--dimensional CR--structures.  In contrast to
the Lagrangian contact structures, there is no distinguished
$G_0$--invariant subset in $\g_{-1}$, so the discussion is similar as
above, but easier, so we skip the details.
\end{example}

\begin{example}\label{4.5}
  Let us finish the paper with the discussion of generalized geodesics
  in the so called x--x--dot geometries (the name comes from the shape
  of the Dynkin diagram with crosses describing the corresponding
  parabolic subgroup in $\frak s\frak l(4,\Bbb C)$). Such structures
  appear as correspondence spaces in classical twistor theory, and
  they are related to the geometric theory of ODE's.

Let us consider the group $G=SL(4,\Bbb R)$ with the parabolic
subgroup $P$ which may be indicated as
$P=\left\{\smatrix{\ast&\ast&\ast&\ast\\0&\ast&\ast&\ast\\
0&0&\ast&\ast\\0&0&\ast&\ast\\}\right\}$.
The following discussion may be also understood as a block--wise
generalization of the discussion of the matrices in the example \ref{4.4}(1)
which we shall call the `x--x' case.
The examples with more `dots' in the Dynkin diagram and just two crosses over
the first two nods on the left will behave quite similarly to the 
x--x--dot case.

The Lie algebra $\g_-$ is described by block matrices of the form
$\g_-=\left\{\smatrix{0&0&0\\x_1&0&0\\X_2&X_1&0}\right\}$, where the blocks
$x_1$, $X_1$ generate the subalgebra $\g_{-1}$ and $X_2$ belongs to
$\g_{-2}$.
The truncated adjoint action of an element 
$\exp\smatrix{0&z_1&Z_2\\0&0&Z_1\\0&0&0}\in P_+$ is given by the formula 
$\smatrix{0&0&0\\x_1&0&0\\X_2&X_1&0}\mapsto
\smatrix{0&0&0\\x_1+Z_1(X_2)&0&0\\X_2&X_1-z_1X_2&0}$.

In accordance with the x--x case, there are two distinguished $G_0$--invariant 
subspaces in $\g_{-1}$ corresponding to the blocks $x_1$ and $X_1$,
respectively. 
The generalized geodesics emanating in the appropriate directions of the
distribution $T^{-1}M$ have got the same properties as above.
In particular, 
curves of this type are determined by a 2--jet but as unparametrized curves
they are given by a direction.
Parametrized geodesics of this type with the common tangent vector form a 
1--dimensional family parametrized by the elements of the form
$\left\{\smatrix{0&z_1&0\\0&0&0\\0&0&0}\right\}$ and 
$\left\{\smatrix{0&0&0\\0&0&Z_1\\0&0&0}\right\}\Big/K$, respectively,
where $K=\left\{\smatrix{0&0&0\\0&0&Z_1\\0&0&0}:Z_1(X_1)=0\right\}$, briefly
written as $K=\{Z_1(X_1)=0\}$.
In the latter case, what really affects on the 2--jet is the value $Z_1(X_1)$
instead of $Z_1$, that is why the quotient appears.

Generalized geodesics with the generic directions in $T^{-1}M$ are determined
by a 3--jet and to any tangent vector there is a 3--dimensional family of
(projectively) parametrized geodesics described by elements of $\p_+/K$, 
where $K=\{z_1=0,Z_1(X_1)=0,Z_2(X_1)=0\}$.

The only contrast with the x--x case appears in the directions not belonging
to $T^{-1}M$. 
The analogy of chains, i.e.\ the curves from $\Cal C_{\g_{-2}}$, does not
exhaust all directions out of the distribution $T^{-1}M$ but only a
4--dimensional `cylinder' 
$\left\{\smatrix{0&0&0\\Z_1(X_2)&0&0\\X_2&-z_1X_2&0}\right\}\subset\g_-$
(at each point) according to the orbit of $\g_{-2}$ with respect to the 
truncated adjoint action of $P$.
Obviously, the complement is formed by all elements of $\g_-$ such that
vectors $X_1$ and $X_2$ are linearly independent; this set is
$G_0$--invariant.
Now, the discussion splits into two branches where the first one follows the 
x--x case, but the second one brings something new.

Let us start with the directions given by chains.
First of all, it is easy to verify that the sets of curves given by the
invariant subsets
$A_1=\left\{\smatrix{0&0&0\\0&0&0\\X_2&aX_2&0}\right\}$ and
$A_2=\left\{\smatrix{0&0&0\\x_1&0&0\\X_2&0&0}\right\}$ are the same and
both of these choices coincide with chains defined by $A=\g_{-2}$. Of
course, all chains depend on 2--jets in one point.
For any tangent vector of this type there is a 1--dimensional family of
parametrized chains, described by the elements of $\g_2/\{Z_2(X_2)=0\}$, 
all parameterizing the same curve.
 
Besides the chains, there is a 3--dimensional family of generalized geodesics
emanating in the same directions as chains from a given point,
defined by the subset $A=\left\{\smatrix{0&0&0\\x_1&0&0\\X_2&aX_2&0}\right\}$.
This family is parametrized by the quotient $\p_+/K$, where
$K=\{z_1=0,Z_1(X_2)=0,Z_2(X_2)=0\}$. 
Curves of this type are also determined by a 2--jet and the admissible 
reparametrizations are affine.

Finally, we fix a tangent vector which does not belong to $T^{-1}M$ and is not 
tangent to a chain.
By analogy to the previous case, there are two disjunct classes of generalized
geodesics emanating in such directions, but having rather different properties
than above.
The first class corresponds to the invariant subset 
$A=\left\{\smatrix{0&0&0\\0&0&0\\X_2&X_1&0}\right\}$, where $X_1$ and $X_2$ 
are supposed to be linearly independent (we assume this in the rest of the 
example). Curves of this type are determined by a 2--jet, they allow projective
reparametrizations, and to the given tangent vector there is a 3--dimensional 
family of parametrized geodesics described by elements of the form
$\left\{\smatrix{0&0&Z_2\\0&0&Z_1\\0&0&0}:Z_1(X_2)=0\right\}$. 
The last distinguished class of curves corresponds to the generic choice of
$A=\left\{\smatrix{0&0&0\\x_1&0&0\\X_2&X_1&0}\right\}$. 
Again, curves of this type are determined by a 2--jet and allow the projective
class of reparametrizations.
The family of parametrized geodesics with the common tangent vector has got
the maximal dimension 5 and it is described by 
all elements of $\p_+$.

\end{example}


\begin{thebibliography}{XX}

\bibitem{BE}
T.N. Bailey, M.G. Eastwood,
Complex paraconformal manifolds: their differential geometry
and twistor theory, Forum Math. {\bf 3} (1991), 61--103.

\bibitem{BE-conf}
T.N. Bailey, M.G. Eastwood,
Conformal circles and parametrizations of curves in conformal manifolds,
Proc. of AMS {\bf 108} (1990), 215--221.

%
%
%
%
%
%
%
\bibitem{CS} A. \v Cap, H. Schichl, Parabolic Geometries and Canonical
  Cartan Connections, Hokkaido Math. J. \textbf{29} No.3 (2000)
  453--505.
 
\bibitem{CS-weyl} A. \v Cap, J. Slov\'ak, Weyl structures for
  parabolic geometries, to appear in Math. Scand. {\bf 93}
  (2003),electronically available as Preprint ESI 801 at www.esi.ac.at

%
\bibitem{CSS4} A. \v Cap, J. Slov\'ak, V. Sou\v cek, 
Bernstein--Gelfand--Gelfand sequences, Annals of Mathematics {\bf 154}
(2001), 97--113.

\bibitem{Cartan-conf} E. Cartan, Les espaces \`a connexion conforme, Ann.
Soc. Pol. Math. {\bf 2} (1923), 171--202.

\bibitem{CheM} S.S. Chern, J. Moser, Real hypersurfaces in complex
manifolds, Acta Math. {\bf 133} (1974), 219--271.

\bibitem{Eas-affine} M.G. Eastwood, V. Ezhov, On affine normal forms and a
classification of homogeneous surfaces in affine three--space, Geometriae
Dedicata {bf 77} (1999), 11--69.

%
\bibitem{Fef} C. Fefferman, Parabolic invariant theory in complex
analysis, Adv. in Math. {\bf 31} (1979), 131--262.

%
%
%
%
%
%
%
%
\bibitem{KMS} I. Kol\'a\v r, P.W. Michor, J. Slov\'ak, {\em Natural
Operations in Differential Geometry\/}, Springer 1993.

\bibitem{Kolar} I. Kol\'a\v r, Higher order torsions of spaces with Cartan
Connection, Cahiers Topologie G\'eom. Diff\'erentielle, 12 (1971), 137--146.

%

%
\bibitem{Koch2} L.K. Koch, Chains, null-chains, and CR geometry,
Trans. Amer. Math. Soc. 338 (1993),  245--261.

\bibitem{Koch3} L.K. Koch, Development and distinguished curves, Preprint, 1993.

%
%
%
\bibitem{Mor} T. Morimoto,  Geometric structures on filtered manifolds,
Hokkaido Math. J. {\bf 22} (1993), 263--347.

\bibitem{Och} T. Ochiai, Geometry associated with semisimple flat
homogeneous spaces, Trans. Amer. Math. Soc. {\bf 152} (1970), 159--193.

\bibitem{SchmSl} G. Schmalz; J. Slov\'ak, The geometry
of hyperbolic and elliptic CR-manifolds of codimension two, Asian J.
Math. {\bf 4} (2000), 565--597.

\bibitem{Sha} R.W. Sharpe, {\em Differential Geometry\/}, Graduate
Texts in Mathematics 166, Springer--Verlag 1997.

%
\bibitem{Tak} M. Takeuchi,
Lagrangian contact structures on projective cotangent bundles,
Osaka J. Math. {\bf 31} (1994), 837--860.

\bibitem{Tan} N. Tanaka,  On the equivalence problem asociated 
with simple graded Lie algebras, Hokkaido Math. J. {\bf 8} (1979), 23--84.

%
%
%
\bibitem{Yam} K. Yamaguchi, Differential systems associated with simple
graded Lie algebras, Advanced Studies in Pure Mathematics {\bf 22} (1993),
413--494.

\bibitem{Zadnik-dissertation} V. \v Z\'adn\'\i k, 
Generalized geodesics, PhD Thesis, Masaryk
University in Brno, 2003.
\end{thebibliography}
\end{document}